\documentclass[12pt]{amsart}
\usepackage{amsmath,amscd,amssymb,amsfonts,graphics,mathrsfs}
\usepackage{xcolor,hyperref}
\hypersetup{colorlinks=true,citecolor=black,linkcolor=black}
\newcommand{\hl}{\hyperlink}
\newcommand{\htt}{\hypertarget}
\newcommand{\h}{\hbox}
\newcommand{\q}{\quad}
\newcommand{\nin}{\noindent}
\newcommand{\bs}{\par\bigskip}
\newcommand{\ms}{\par\medskip}
\newcommand{\sk}{\par\smallskip}
\newcommand{\bsn}{\par\bigskip\noindent}
\newcommand{\msn}{\par\medskip\noindent}

\newcommand{\ges}{\geqslant}
\newcommand{\les}{\leqslant}
\newcommand{\1}{\hskip1pt}

\newcommand{\mcup}{\hbox{$\bigcup$}}
\newcommand{\mopl}{\hbox{$\bigoplus$}}
\newcommand{\msum}{\hbox{$\sum$}}

\newcommand{\B}{{\mathscr B}}
\newcommand{\D}{{\mathscr D}}
\newcommand{\Fc}{{\mathscr F}}
\newcommand{\Lc}{{\mathscr L}}
\newcommand{\OO}{{\mathscr O}}
\newcommand{\Bf}{{\mathscr B}_{\!f}}
\newcommand{\Bt}{\widetilde{\mathscr B}}
\newcommand{\Btf}{\widetilde{\mathscr B}_{\!f}}

\newcommand{\Mtf}{{}\hskip9pt\widetilde{\hskip-9pt\mathscr M}{}_{\hskip-4pt f}}
\newcommand{\Vf}{{\mathscr V}_{\!f}}
\newcommand{\VX}{{\mathscr V}_{\!X}}
\newcommand{\Vt}{\widetilde{V}}

\newcommand{\Xt}{\widetilde{X}}
\newcommand{\Yt}{\widetilde{Y}}
\newcommand{\alt}{\widetilde{\alpha}}

\newcommand{\PP}{{\mathbb P}}
\newcommand{\Q}{{\mathbb Q}}
\newcommand{\C}{{\mathbb C}}
\newcommand{\Fb}{{\mathbb F}}
\newcommand{\N}{{\mathbb N}}
\newcommand{\DD}{{\mathbf D}}

\newcommand{\Z}{{\mathbb Z}}
\newcommand{\BM}{{\rm BM}}
\newcommand{\Gr}{{\rm Gr}}
\newcommand{\DR}{{\rm DR}}
\newcommand{\Sing}{{\rm Sing}}
\newcommand{\ual}{^{(\al)}}
\newcommand{\al}{\alpha}
\newcommand{\be}{\beta}

\newcommand{\la}{\lambda}

\newcommand{\Om}{\Omega}
\newcommand{\si}{\sigma}

\newcommand{\dd}{\partial}

\newcommand{\eq}{\,{=}\,}
\newcommand{\gess}{\,{\ges}\,}
\newcommand{\less}{\,{\les}\,}
\newcommand{\sgt}{\,{>}\,}
\newcommand{\slt}{\,{<}\,}
\newcommand{\nes}{\,{\ne}\,}

\newcommand{\mi}{\1{-}\1}

\newcommand{\pl}{\1{+}\1}

\newcommand{\bl}{\bigl}
\newcommand{\br}{\bigr}
\newcommand{\sst}{\,{\subset}\,}
\newcommand{\stm}{\,{\setminus}\,}
\newcommand{\ins}{\,{\in}\,}
\newcommand{\tos}{\,{\to}\,}
\newcommand{\defs}{\,{:=}\,}
\newcommand{\ssc}{\,\raise.15ex\hbox{${\scriptstyle\circ}$}\,}
\newcommand{\ssb}{\raise.15ex\h{${\scriptscriptstyle\bullet}$}}
\newcommand{\into}{\hookrightarrow}

\newcommand{\simto}{\,\,\rlap{\hskip1.3mm\raise1.4mm\hbox{$\sim$}}\hbox{$\longrightarrow$}\,\,}
\newcommand{\dX}{\delta\!_X}
\begin{document}
\h{}\bs
\centerline{\large Hirzebruch-Milnor classes of hypersurfaces with nontrivial normal bundles}
\sk
\centerline{\large and applications to higher du~Bois and rational singularities}
\bs
\centerline{Lauren\c{t}iu G. Maxim, Morihiko Saito, and Ruijie Yang}
\bsn\ms
\vbox{\nin\narrower\smaller
{\bf Abstract.} We extend the Hirzebruch-Milnor class of a hypersurface $X$ to the case where the normal bundle is nontrivial and $X$ cannot be defined by a global function, using the associated line bundle and the graded quotients of the monodromy filtration. The earlier definition requiring a global defining function of $X$ can be applied rarely to projective hypersurfaces with non-isolated singularities. Indeed, it is surprisingly  difficult to get a one-parameter smoothing with total space smooth without destroying the singularities by blowing-ups (except certain quite special cases). As an application, assuming the singular locus is a projective variety, we show that the minimal exponent of a hypersurface can be captured by the spectral Hirzebruch-Milnor class, and higher du~Bois and rational singularities of a hypersurface are detectable by the unnormalized Hirzebruch-Milnor class. Here the unnormalized class can be replaced by the normalized one in the higher du~Bois case, but for the higher rational case, we must use also the decomposition of the Hirzebruch-Milnor class by the action of the semisimple part of the monodromy (which is equivalent to the spectral Hirzebruch-Milnor class). We cannot extend these arguments to the non-projective compact case by Hironaka's example.}
\ms\bs
\centerline{\bf Introduction}
\bsn
Let $X$ be a reduced complex algebraic variety. For $k\ins\N$, the notions of $k$-du~Bois and $k$-rational singularities are introduced in \cite{JKSY3}, \cite{FL} generalizing the classical notions for $k\eq0$. These can be characterized for hypersurfaces respectively by the conditions
\htt{1}{}
$$\alt_X\gess k{+}1\q\q\h{and}\q\q\alt_X\sgt k{+}1,
\leqno(1)$$
\par\nin extending the case $k\eq0$ in \cite{rat}, \cite{mos}. Here the {\it minimal exponent\1} $\alt_X$ is defined to be the minimum of $\alt_f$ with $f$ running over local defining functions of $X$, and $\alt_f$ denotes up to sign the maximal root of the reduced Bernstein-Sato polynomial $b_f(s)/(s{+}1)$, see \cite{MOPW}, \cite{JKSY3} and \cite[Appendix]{FL}, \cite{MP} respectively. (Some argument in \cite{JKSY3} is unnecessarily too complicated, see \cite[Appendix]{FL}.) As is pointed out in \cite{MY} the above characterizations imply immediately the Thom-Sebastiani type assertions by using \cite[(3.2.3)]{MSS2} or \cite[Thm.\,0.8]{mic} in the case one function is weighted homogeneous.
\sk
For a hypersurface $X$ in a smooth variety $Y$, let
$$M_y(X)\ins H^{\BM}_{2\ssb}(\Sing\,X,\Q)[y]$$
\par\nin be the {\it normalized Hirzebruch-Milnor class,} see \cite{MSS1}. In case the normal bundle is nontrivial, this can be defined by using the line bundle on $Y$ associated with $X$ (instead of the graph embedding, see for instance \cite[2.2]{mhm}) and the graded quotients of the {\it monodromy filtration,} see \hl{1.4}{1.4} below. The earlier definition requiring a {\it global\1} defining function of $X$ (see \cite{MY}) can be applied {\it rarely\1} to projective hypersurfaces with {\it non-isolated\1} singularities. Indeed, we can construct a good smoothing, that is, with total space smooth, without destroying the singularities by blowing-ups in the case where the singular locus is one-dimensional and the general transversal slices are ordinary $m$-ple points, but it is extremely difficult to relax these conditions, see \cite{ex}. (Note that the total space of a one-parameter smoothing given by a linear system is always singular, since $V_g\,{\cap}\,\Sing\,V_f\sst\Sing\,V_{f+sg}$ with $V_f\defs f^{-1}(0)$, see Remark~\hl{R1.4c}{1.4c} below, and we get {\it reducible\1} varieties having {\it irrational\1} and {\it non-$1$-du~Bois\1} singularities by taking the {\it total\1} transform if we desingularize the total space by blowing-ups with smooth centers contained in the singular locus. Note also that the spectral Hirzebruch-Milnor class of an {\it isolated\1} hypersurface singularity is equivalent to the spectrum.) So this paper supersedes the preceding one \cite{MY} in this point. This change of formulation seems to be a rather {\it conceptual\1} one, since one cannot replace it by a \v Cech type argument using the vanishing cycle functor, see Remark~\hl{R1.4b}{1.4b} below.
\sk
We denote by $M_y(X)^{(p)}\ins H^{\BM}_{2\ssb}(\Sing\,X,\Q)$ the coefficient of $y^p$ in $M_y(X)$ so that
$$M_y(X)\eq\msum_{p\in\N}\,M_y(X)^{(p)}\1y^p.$$
\par\nin Let $M^{\rm un}_y(X)$ be the {\it unnormalized Hirzebruch-Milnor class\1} of $X$. By definition we have the relation
\htt{2}{}
$$M^{\rm un}_y(X)_{2k}\eq(1{+}\1y)^kM_y(X)_{2k},
\leqno(2)$$
\par\nin where $M_y^{\rm un}(X)_{2k}\ins H^{\BM}_{2k}(\Sing\,X,\Q)[y]$ denotes the degree $2k$ part of $M_y^{\rm un}(X)$, and similarly for $M_y(X)_{2k}$. We can define $M_y(X)^{(p)}_{2k}$ similarly. The equality (\hl{2}{2}) shows that $M^{\rm un}_y(X)_{2k}$ is divisible by $(1{+}y)^k$, see \cite[Prop.\,5.21]{Sch}. It is easy to verify this for the highest degree non-zero part using \hl{2.2}{2.2} below, since the ranks of components of the de Rham complex of a smooth variety are given by {\it binomial coefficients.}
\sk
There is the decomposition by the action of the semisimple part of the monodromy
\htt{3}{}
$$M_y(X)\eq\msum_{|\la|=1}\1M_y^{\{\la\}}(X)\q\h{in}\,\,\,\,H^{\BM}_{2\ssb}(\Sing\,X,\Q)[y],
\leqno(3)$$
\par\nin (similarly for $M_y^{\rm un}(X)$). This is essentially equivalent to considering the {\it spectral Hirzebruch-Milnor class,} see \cite{MSS3} (and (\hl{8}{8}) below).
\sk
In this paper we prove the following.
\par\htt{T1}{}\msn
{\bf Theorem~1.} {\it Let $X$ be a hypersurface of a smooth variety $Y$, and $k\ins\N$. If $X$ has only $k$-du~Bois singularities, then
\htt{4}{}
$$M_y(X)^{(p)}\eq0\q(\forall\,p\less k),
\leqno(4)$$
\par\nin \vskip-7mm
\htt{5}{}
$$M_y^{\rm un}(X)^{(p)}\eq0\q(\forall\,p\less k),
\leqno(5)$$
\par\nin where the last two properties are equivalent to each other by {\rm(\hl{2}{2})}. The converse holds if $\,\Sing\,X$ is a projective variety.}
\par\htt{T2}{}\msn
{\bf Theorem~2.} {\it Let $X$ be a hypersurface of a smooth variety $Y$, and $k\ins\N$. If $X$ has only $k$-rational singularities, then
\htt{6}{}
$$M_y(X)^{(p)}\eq0\q(\forall\,p\less k),\q M_y^{\{1\}}(X)^{(k+1)}\eq0,
\leqno(6)$$
\par\nin \vskip-7mm
\htt{7}{}
$$M_y^{\rm un}(X)^{(p)}\eq0\q(\forall\,p\gess d_X{-}k).
\leqno(7)$$
\par\nin The converse holds assuming either {\rm(\hl{6}{6})} or {\rm(\hl{7}{7})}
if $\,\Sing\,X$ is a projective variety.}
\ms
It does not seem quite clear whether (\hl{6}{6}) and (\hl{7}{7}) are equivalent to each other without assuming the projectivity of $\Sing\,X$.
\par\htt{R1}{}\msn
{\bf Remark~1.} We cannot replace {\it projectivity\1} by {\it compactness\1} in the assumption of the converse assertions because of Hironaka's example \cite{Hi}, \cite[Appendix B, Example 3.4.1]{Ha}, see Proposition~\hl{P2}{2} and \hl{2.3}{2.3} below. Note that the {\it positivity\1} of the Todd class transformation in the projective case (see \hl{2.2}{2.2} below) is crucial to the proof of the converse assertions in the main theorems. We cannot omit the {\it compactness\1} assumption, since the Chow groups of affine varieties are usually small (for instance, $\C^n$).
\par\htt{R2}{}\msn
{\bf Remark~2.} Some assertions similar to Theorems~\hl{T1}{1} and \hl{T2}{2} for {\it unnormalized\1} classes have been given initially in \cite{MY}, where $X$ was assumed to be {\it globally\1} defined by a function $f$ and the {\it unnormalized\1} version of (\hl{6}{6}) together with (\hl{5}{5}) was treated. Note that duality is not needed for (\hl{7}{7}).
\par\htt{R3}{}\msn
{\bf Remark~3.} The main theorems are quite {\it theoretical,} and are not necessarily good for {\it explicit computations,} since it is usually easier to calculate locally the {\it minimal exponents.} This is the reason for which we use the word ``detect". It is important that Theorems~\hl{T1}{1} and \hl{T2}{2} are {\it directly applicable to any hypersurfaces of smooth projective varieties\1} without taking any blow-ups nor assuming a condition on $k$ as in \cite[Thm.\,1.3]{ex} which is written to show the ``limitation" of the old definition of Hirzebruch-Milnor class. Here a calculation for the detection of non-$k$-du~Bois (or rational) singularities is {\it essentially local\1} (so the divisor is principal) since we can restrict {\it inductively\1} to an open subset of the {\it possibly\1} non-$k$-du~Bois (or rational) locus, using a stratification together with the {\it positivity\1} of the Todd class transformation in the {\it projective\1} case, see Proposition~\hl{P1}{1} and \hl{2.2}{2.2} below.
\sk
For the proof of the converse assertions in Theorems~\hl{T1}{1} and \hl{T2}{2}, we need for instance the following.
\par\htt{P1}{}\msn
{\bf Proposition~1.} {\it If $Z$ is a projective variety and $\Fc$ is a coherent sheaf on $Z$ whose support has dimension $k$, then $td_*[\Fc]_{2k}\ins H_{2k}(Z,\Q)$ does not vanish.}
\sk
This can be reduced to the case $Z\eq\PP^n$ by using the compatibility of the Todd class transformation $td_*$ with the pushforward under proper morphisms, and follows from the theory of {\it topological filtration,} where the assertion is reduced to the {\it positivity\1} of the degrees of subvarieties of $\PP^n$ (which is defined by taking the intersection with a sufficiently general linear subspace of complementary dimension), see for instance \cite{SGA6}, \cite[Cor.\,18.3.2]{Fu}, and \cite[1.3 (or 1.6 in the preprint version)]{MSS3} (and also \hl{2.2}{2.2} below).
\sk
One cannot detect an element of the Grothendieck group of coherent sheaves by applying the Todd class transformation without assuming that the coefficients of coherent sheaves with highest-dimensional supports have the {\it same sign.} The last condition is satisfied in the case of coherent sheaves and also for the dual of coherent sheaves ({\it locally}). So an {\it inductive\1} argument is needed, see \hl{2.1}{2.1} below for details.
\sk
In case the singular locus is a projective variety, the {\it minimal exponent\1} $\alt_X$ of a hypersurface $X\sst Y$ can be detected by using the  {\it normalized spectral Hirzebruch-Milnor class\1}
\htt{8}{}
$$M_t^{\rm sp}(X)\defs\msum_{|\la|=1,\,p\in\Z}\,(-1)^pM_y^{\{\la\}}(X)^{(p)}t^{\1p+\ell(\la)}\in H^{\rm BM}_{2\ssb}({\rm Sing}\,X,\Q)[t^{1/e}],
\leqno(8)$$
\par\nin (similarly for the {\it unnormalized\1} one with $M_y$ replaced by $M_y^{\rm un}$), see \cite{MSS3}. Here $e$ is a positive integer divisible by the order of the semisimple part of the monodromy on the vanishing cohomologies, and $\ell(\la)\ins[0,1)$ is defined by $e^{2\pi i\ell(\la)}\eq\la$ for $\la\ins\C^*_1$, setting $\C^*_1\defs\{\la\ins\C^*\,|\,|\la|\eq1\}$. We have the following.
\par\htt{T3}{}\msn
{\bf Theorem~3.} {\it Let $X$ be a hypersurface of a smooth variety $Y$, and $\al\ins\Q$. If $\al\slt\alt_X$, then
\htt{9}{}
$$M_y^{\{\la\}}(X)^{(p)}\eq0\q\h{for any $\,p,\la\,$ with $\,p\pl\ell(\la)\less\al$}.
\leqno(9)$$
\par\nin The converse holds in the case ${\rm Sing}\,X$ is a projective variety. Moreover these assertions remain true with normalized class $M_y(X)$ replaced by the unnormalized one $M_y^{\rm un}(X)$.}
\ms
The last assertion follows from the divisibility in (\hl{2}{2}). Theorems~\hl{T1}{1} and \hl{T2}{2} except (\hl{7}{7}) are corollaries of Theorem~\hl{T3}{3}.
\sk
We define more precisely the $\lambda$-minimal exponents by
\htt{10}{}
$$\alt_X^{\{\la\}}\defs\min\bl\{\al\ins\Q\mid\Gr_{\Vt}^{\al}\OO_Y\nes0,\,e^{-2\pi i\al}\eq\la\br\}\q\h{for}\,\,\,\la\ins\C^*_1,
\leqno(10)$$
\par\nin in a compatible way with (\hl{1.1.3}{1.1.3}) below. These numbers can be captured by the spectral Hirzebruch class, since we have the following.
\par\htt{T4}{}\msn
{\bf Theorem~4.} {\it Let $X$ be a hypersurface of a smooth variety $Y$, and $\al\ins\Q$. Set $\la=e^{2\pi i\al}$. If $\al\slt\alt_X^{\{\la^{-1}\}}$, then
\htt{11}{}
$$M_y^{\{\la\}}(X)^{(p)}\eq0\q\h{for any $\,p\in\Z\,$ with $\,p\pl\ell(\la)\less\al$}.
\leqno(11)$$
\par\nin The converse holds in the case ${\rm Sing}\,X$ is a projective variety. Moreover these assertions remain true with normalized class $M_y(X)$ replaced by the unnormalized one $M_y^{\rm un}(X)$.}
\ms
Theorem~\hl{T3}{3} is a corollary of Theorem~\hl{T4}{4}, since $\alt_X\eq\min\bl\{\alt_X^{\{\la\}}\br\}$. The $\la$-minimal exponents $\alt_f^{\{\la\}}$ cannot be determined by the Bernstein-Sato polynomial because of the shift of roots under $\mu$-constant deformations, see for instance \cite{sem}. Note that the spectral Hirzebruch-Milnor class in \cite{MSS3} corresponds to the spectrum as in \cite{GLM}, which is denoted by ${\rm Sp}'_f(t)$ in \cite{ste} and is the ``dual" of ${\rm Sp}_f(t)$, more precisely, ${\rm Sp}'_f(t)\eq{\rm Sp}_f(t^{-1})t^n$. This is closely related to the self-duality used in an essential way for the proofs of Theorem~\hl{T4}{4}, see (\hl{1.1.4}{1.1.4}--\hl{1.1.5}{5}) below.
\sk
Concerning the converse assertions in Theorems~\hl{T1}{1}--\hl{T3}{3}, there is a slightly stronger version as below. Let $\dX$ be the dimension of the {\it minimal exponent locus\1} $\{x\ins X\mid\alt_{X,x}\eq\alt_X\}$, where $\alt_{X,x}$ is the minimal exponent of $X$ at $x$ defined {\it up to sign} by the maximal root of the reduced local Bernstein-Sato polynomial $b_{f,x}(s)/(s{+}1)$. We have the following.
\par\htt{T5}{}\msn
{\bf Theorem~5.} {\it Assume $X$ is a projective hypersurface. Then the image of $M_y^{\{\la\}}(X)^{(p)}_{2\dX}$ in $H_{2\dX}(X,\Q)$ does not vanish for $p\ins\Z$, $\la\ins\C^*_1$ with $p\pl\ell(\la)=\alt_X$. Similarly the converse assertions in Theorems~{\rm\hl{T1}{1}} and {\rm\hl{T2}{2}} hold by replacing $M_y(X)^{(p)}$,  $M^{\rm un}_y(X)^{(p)}$ with the images of their degree $2\dX$ part in $H_{2\dX}(X,\Q)$.}
\ms
These stronger assertions, however, fail in the {\it non-projective\1} compact case by using a variant of Hironaka's example \cite{Hi} as follows (see \hl{2.3}{2.3} below).
\par\htt{P2}{}\msn
{\bf Proposition~2.} {\it There is a non-projective compact hypersurface $X$ such that the assertions of Theorem~{\rm\hl{T5}{5}} do not hold.}
\ms
As for the relation between the spectral Hirzebruch-Milnor class and the Hodge ideals in the sense of Mustata and Popa, this seems quite nontrivial since their Hodge ideals $I_k(\alpha X)$ are not necessarily weakly decreasing for the index $\alpha$ (see \cite{JKSY1}), so the notion of jumping coefficients cannot be defined easily for them.
\sk
At the end we would like to clarify the relation of the present paper with the earlier one \cite{MY} by the first and third named authors. In the latter the authors considered only the case of hypersurfaces in smooth varieties defined by global functions. Soon after the paper appeared on the web, the second named author generalized (spectral) Hirzebruch-Milnor classes to arbitrary hypersurfaces in smooth varieties, and extended some assertions in \cite{MY} to this setting. The current paper is intended to bring the most general results to light; as such, it supersedes \cite{MY}.
\sk
In Section 1 we explain the Hirzebruch-Milnor classes of hypersurfaces with nontrivial normal bundles after reviewing some basics in Hodge module theory. In Section 2 we prove the main theorems using the characterizations of higher du~Bois and rational singularities via the minimal exponents. Here the topological filtration on the Grothendieck group of coherent sheaves plays an essential role for the converse assertions, and its short account is also given.
\msn
{\bf Acknowledgement.} The first named author is partially supported by the Simons Foundation (Collaboration Grant \#567077), and by the Romanian Ministry of National Education (CNCS-UEFISCDI grant PN-III-P4-ID-PCE-2020-0029).
The second named author was partially supported by JSPS Kakenhi 15K04816.
\bs
\centerline{Contents}
\ms
\par\hl{S1}{1. Hirzebruch-Milnor classes of hypersurfaces}\hfill 5
\par\q\hl{1.1}{1.1.~Vanishing cycle filtered $\D$-modules}\hfill 5
\par\q\hl{1.2}{1.2.~Algebraic microlocalization}\hfill 6
\par\q\hl{1.3}{1.3.~Microlocal $V$-filtration}\hfill 6
\par\q\hl{1.4}{1.4.~Hirzebruch-Milnor classes}\hfill 7
\par\hl{S2}{2. Proofs of the main theorems}\hfill 7
\par\q\hl{2.1}{2.1.~Proofs of Theorems}~\hl{T1}{1}, \hl{T2}{2}, \hl{T4}{4} and \hl{T5}{5}\hfill 7
\par\q\hl{2.2}{2.2.~Topological filtration}\hfill 8
\par\q\hl{2.3}{2.3.~A variant of Hironaka's example}\hfill 9
\par\bs\bs\htt{S1}{}
\vbox{\centerline{\bf 1. Hirzebruch-Milnor classes of hypersurfaces}
\bsn
In this section we explain the Hirzebruch-Milnor classes of hypersurfaces with nontrivial normal bundles after reviewing some basics in Hodge module theory.}
\par\htt{1.1}{}\msn
{\bf 1.1.~Vanishing cycle filtered $\D$-modules.} Let $Y$ be a smooth complex algebraic variety with $X\sst Y$ a reduced hypersurface. Let $Z\tos Y$ be the line bundle corresponding to the invertible sheaf $\OO_Y(X)$. Here $Y$ can be identified with the zero-section and also with the canonical section $s_{X,Y}$ corresponding to $1\in\OO_Y(X)$, and the inclusion morphism of the latter is denoted by $i_{Y,Z}:Y\into Z$. This replaces the graph embedding $i_f$ in the case $X$ is defined by $f$. Note that the zero-locus of $s_{X,Y}$ is $X$.
\sk
Let $\B_{Y,Z}$ be the direct image of $\OO_Y$ by $i_{Y,Z}$ as a left $\D$-module. Choosing a local defining equation $f$ of $X$, we have locally an isomorphism
\htt{1.1.1}{}
$$\B_{Y,Z}|_U\eq\Bf\,({=}\,\OO_U[\dd_t]).
\leqno(1.1.1)$$
\par\nin Here $\Bf$ is the direct image of $\OO_U$ by the graph embedding $i_f:U\into U{\times}\C$ with $t$ the coordinate of $\C$, and $f$ is a defining function of $X$ on an open subvariety $U\sst Y$.
\sk
We define the Hodge filtration $F$ by the order of $\dd_t$ so that there are locally isomorphisms
\htt{1.1.2}{}
$$\Gr^F_0\B_{Y,Z}|_U\eq\Gr^F_0\Bf\eq\OO_U.
\leqno(1.1.2)$$
\par\nin Note that this filtration is shifted by 1 compared with the filtration defined by the direct image under the closed embedding $i_{Y,Z}$. This shift is indispensable for the compatibility of the Hodge filtration on the de Rham complex with the direct image.
\sk
Let $V$ be the filtration of Kashiwara \cite{Ka} and Malgrange \cite{Ma} indexed by $\Q$ on $\B_{Y,Z}$ along the zero-section (where it is not necessary to assume that the line bundle is trivial). Here $\theta{-}\al{+}1$ is nilpotent on $\Gr_V^{\al}$ with $\theta$ the Euler field corresponding to the $\C^*$-action on the fibers of the line bundle. Recall that $\B_{Y,Z}$ is locally identified with $\Bf$.
\par\htt{D1.1}{}\msn
{\bf Definition~1.1.} We define the {\it vanishing cycle filtered $\D$-module\1} $(\VX,F)\eq\mopl_{|\la|=1}\,(\VX^{\{\la\}},F)$ as follows:
\htt{1.1.3}{}
$$\aligned(\VX^{\{\la\}},F)&\defs\Gr_V^{\al}(\B_{Y,Z},F[1])\q(\la\eq e^{-2\pi i\al},\,\al\ins(0,1)),\\(\VX^{\{1\}},F)&\defs\Gr_V^{0}(\B_{Y,Z},F).\endaligned
\leqno(1.1.3)$$
\par\nin Here $\VX^{\{\la\}}\eq0$ unless $\la$ is a root of unity, and $(F[m])_p\eq F_{p-m}$ in a compatible way with $F_p\eq F^{-p}$ ($p,m\ins\Z$). We denote $\VX|_U$ by $\Vf$ when a local defining function $f$ of $X$ is chosen on an open subvariety $U\sst Y$.
\sk
By definition the $\VX^{\{\la\}}$ are filtered $\Gr_V^{0}\D_Z$-modules, but there are no canonical structures of $\D_Y$-modules, and $\Gr_V^0\D_Z$ is {\it locally\1} isomorphic to $\D_Y[\theta]$ (choosing $f$) with $\theta$ the Euler field corresponding to the $\C^*$-action on the fiber of line bundle. If the line bundle is trivialized, $V^0\D_Z$ is generated by $\D_Y$, $\OO_Z$, and $\theta$ as a ring, and $V^k\D_Z\eq t^kV^0\D_Z$ for $k\in\N$. Note that $\Vf^{\{\la\}}$ has a structure of a filtered $\D_U$-module (choosing $f$).
\par\htt{R1.1}{}\msn
{\bf Remark~1.1.} The above construction is closely related to {\it Verdier specialization\1} \cite{Ve}, see for instance \cite[1.3]{BMS}. It seems that \cite{Ka} is influenced by it, looking at the reference [1] in it.
\ms
The following self-duality is proved in \cite{dual}, and is used in Hodge module theory in an essential way (see also \cite[(2.1.4)]{JKSY2}, \cite{KLS}):
\par\htt{P1.1}{}\msn
{\bf Proposition~1.1.} {\it There are self-duality isomorphisms of filtered left $\D_U$-modules}
\htt{1.1.4}{}
$$\aligned\DD\bl(\Vf^{\{\la\}},F\br)&\eq\bl(\Vf^{\{\la^{-1}\}},F[d_X]\br)\q(\la\nes1),\\ \DD\bl(\Vf^{\{1\}},F\br)&\eq\bl(\Vf^{\{1\}},F[d_Y]\br).\endaligned
\leqno(1.1.4)$$
\par\nin \ms
Here $d_X\defs\dim X\eq d_Y{-}1$, and $\DD$ is the dual functor for filtered left $\D$-modules (which is defined in a compatible way with the one for filtered right $\D$-modules, see for instance \cite{ypg}, \cite{JKSY2}, \cite{JKSY3}). The dual functor $\DD$ is compatible with the de Rham functor $\DR_Y$ and also with $\Gr_F^{\ssb}$ (where $F^p\eq F_{-p}$), see \cite[Ch.\,2]{mhp}. This is used for instance in the calculation of {\it generating levels,\1} see \cite[Remark\,(ii) before Proposition~1.4]{mos}. We then get the following.
\par\htt{C1.1}{}\msn
{\bf Corollary~1.1.} {\it There are self-duality isomorphisms in $D^b_{\rm coh}(\OO_U)\1{:}$}
\htt{1.1.5}{}
$$\aligned\DD\bl(\Gr^p_F\DR_U\bl(\Vf^{\{\la\}}\br)\br)&\eq\Gr^{d_X-p}\DR_U(\Vf^{\{\la^{-1}\}}\br)\q(\la\nes1),\\ \DD\bl(\Gr^p_F\DR_U\bl(\Vf^{\{1\}}\br)\br)&\eq\Gr^{d_Y-p}\DR_U\bl(\Vf^{\{1\}}\br).\endaligned
\leqno(1.1.5)$$
\par\nin \par\htt{1.2}{}\msn
{\bf 1.2.~Algebraic microlocalization.} Let $\Bt_f$ be the {\it algebraic partial microlocalization\1} of $\Bf$ by $\dd_t$, that is,
\htt{1.2.1}{}
$$\Bt_f\eq\OO_U[\dd_t,\dd_t^{-1}],
\leqno(1.2.1)$$
\par\nin see for instance \cite{mic}. The Hodge filtration $F$ is defined by the order of $\dd_t$.
\sk
The following is proved in \cite[2.1--2]{mic} (see also \cite[1.1]{JKSY2}):
\par\htt{P1.2}{}\msn
{\bf Proposition~1.2.} {\it There is the filtration $V$ on $\Bt_f$ such that}
\htt{1.2.2}{}
$$\Gr_V^{\al}(\Bf,F)\simto\Gr_V^{\al}(\Btf,F)\q(\forall\,\al<1),
\leqno(1.2.2)$$
\par\nin \vskip-6mm
\htt{1.2.3}{}
$$\dd_t^j:F_pV^{\al}\Btf\simto F_{p+j}V^{\al-j}\Btf\q(\forall\,j,p\ins\Z,\,\al\ins\Q).
\leqno(1.2.3)$$
\par\nin \par\htt{D1.2}{}\msn
{\bf Definition~1.2.} Set
\htt{1.2.4}{}
$$(\Mtf\ual,F)\defs\Gr_V^{\al}(\Btf,F)\q(\al\ins(0,1]).
\leqno(1.2.4)$$
\par\nin \sk
By Proposition~\hl{P1.2}{1.2}, we get the following.
\par\htt{C1.2}{}\msn
{\bf Corollary~1.2.} {\it There are isomorphisms of filtered left $\D_U$-modules}
\htt{1.2.5}{}
$$(\Vf^{\{\la\}},F)\eq(\Mtf\ual,F[1])\q(\la\eq e^{-2\pi i\al},\,\al\ins(0,1]).
\leqno(1.2.5)$$
\par\nin \par\htt{1.3}{}\msn
{\bf 1.3.~Microlocal $V$-filtration.} In the notation of \hl{1.2}{1.2}, the {\it microlocal $V\!$-filtration\1} $\Vt$ on $\OO_U$ is defined by
\htt{1.3.1}{}
$$(\OO_U,\Vt)\defs\Gr^F_0(\Btf,V).
\leqno(1.3.1)$$
\par\nin \sk
From Proposition~\hl{P1.2}{1.2}, we can deduce the following.
\par\htt{C1.3}{}\msn
{\bf Corollary~1.3.} {\it There are isomorphisms of $\OO_U$-modules}
\htt{1.3.2}{}
$$\Gr_{\Vt}^{\al+p}\OO_U\eq\Gr^F_p\Mtf\ual\q(\forall\,p\ins\Z,\,\al\ins(0,1]).
\leqno(1.3.2)$$
\par\nin \sk
The following is proved in \cite[(1.3.8)]{hi}:
\par\htt{P1.3}{}\msn
{\bf Proposition~1.3.} {\it We have the equality}
\htt{1.3.3}{}
$$\alt_f\eq\min\{\al\in\Q\mid\Gr_{\Vt}^{\al}\OO_U\nes 0\}.
\leqno(1.3.3)$$
\par\nin \sk
Recall that the {\it minimal exponent\1} $\alt_f$ is the absolute value of the maximal root of the reduced Bernstein-Sato polynomial $b_f(s)/(s{+}1)$.
\par\htt{R1.3}{}\msn
{\bf Remark~1.3.} The microlocal $V$-filtration does not depend on the choice of a local defining function $f$. Indeed, if $f'\eq uf$ with $u$ a unit, this $u$ defines the isomorphism of the trivialized line bundles compatible with the graph embeddings by $t'\eq ut$, where $t,t'$ are the coordinates of the trivialized line bundles by $f$ and $f'$ respectively. We have the equality $\dd_t\eq u\dd_{t'}$, since $t\dd_t\eq t'\dd_{t'}$ considering the canonical $\C^*$-action on the line bundles.
\par\htt{1.4}{}\msn
{\bf 1.4.~Hirzebruch-Milnor classes.} In the notation of \hl{1.1}{1.1}, we define $M_y^{{\rm un}\{\la\}}(X)^{(p)}$ in the introduction as follows: Set
\htt{1.4.1}{}
$$M_y^{{\rm un}\{\la\}}(X)^{(p)}\defs(-1)^{p+d_X}\1\msum_{k\in\Z}\,td_*\bl[\Gr_F^p\DR_Y\bl(\Gr^W_k\VX^{\{\la\}}\br)\br],
\leqno(1.4.1)$$
\par\nin where
\vskip-7mm
$$td_*:K_0(\Sing\,X)\to H^{\rm BM}_{2\ssb}(\Sing\,X,\Q)$$
\par\nin is the Todd class transformation in \cite{BFM}, and $W$ is the monodromy filtration associated with the nilpotent action of $\theta{-}\al{+}1$ on $\VX^{\{\la\}}$, where $\al\ins(0,1]$, $\la\eq e^{-2\pi i\al}$. This definition coincides with the one in \cite{MSS1} in the case $X$ is defined {\it globally\1} by a function. Recall that $\theta$ is the Euler field of the line bundle, and $\VX^{\{\la\}}$ has only a structure of a filtered $\Gr_V^0\D_Z$-module.
\sk
The well-definedness of the right-hand-side of {\rm(\hl{1.4.1}{1.4.1})} follows from Proposition~\hl{P1.4}{1.4} just below using for instance the definition of the filtered de Rham complex associated with an integrable connection.
\par\htt{P1.4}{}\msn
{\bf Proposition~1.4.} {\it The $\Gr^F_{\ssb}\Gr^W_k\VX^{\{\la\}}$ have globally well-defined structures of $\Gr^F_{\ssb}\!\1\D_Y\!\1$-modules.}
\msn
{\it Proof.} We first see that the ambiguity of a lifting of a vector field on $Y$ to $V^0\D_Z$ is given by $\OO_{\!Z}\1\theta$. By the definition of the monodromy filtration, the action of $[\theta]$ on $\Gr^F_{\ssb}\Gr^W_k\VX^{\{\la\}}$ is induced by multiplication by $\al\mi1$. This can be neglected since $[\theta]\ins\Gr^F_1\D_Z$. So the assertion follows. This finishes the proof of Proposition~\hl{P1.4}{1.4}.
\par\htt{R1.4a}{}\msn
{\bf Remark~1.4a.} Let $I(\be)$ be the ideal of $\Gr_V^0\D_Z$ generated by $\theta{-}\be$ which is in the center of $\Gr_V^0\D_Z$. The quotient ring $\Gr_V^0\D_Z/I(\be)$ does not seem to be isomorphic to $\D_Y$ in general unless $\be\eq0$.
\par\htt{R1.4b}{}\msn
{\bf Remark~1.4b.} One cannot replace the above construction by an \v Cech type argument using the vanishing cycle functor. Indeed, for instance, the information of the degree of a line bundle in the 1-dimensional singular locus case is lost by taking a \v Cech complex if the differential is forgotten in the Grothendieck group.
\par\htt{R1.4c}{}\msn
{\bf Remark~1.4c.} Let $Y$ be a smooth complex projective variety, and $\Lc$ be an ample line bundle. Let $V_f$ be the zero-locus of $f\ins\Gamma(Y,\Lc)$. If $V_f$ has non-isolated singularities, then $V_{f+sg}\sst Y{\times}\C$ is singular for any $g\ins\Gamma(Y,\Lc)$, where $s$ is the coordinate of $\C$. Indeed,
\htt{1.4.2}{}
$$V_g\cap\Sing\,V_f=\Sing\,V_{f+sg}\cap\{s\eq0\},
\leqno(1.4.2)$$
\par\nin and the left-hand side is non-empty, since $\Lc$ is ample.
\par\bs\bs\htt{S2}{}
\vbox{\centerline{\bf 2. Proofs of the main theorems}
\bsn
In this section we prove the main theorems using the characterizations of higher du~Bois and rational singularities via the minimal exponents. Here the topological filtration on the Grothendieck group of coherent sheaves plays an essential role for the converse assertions, and its short account is also given.}
\par\htt{2.1}{}\msn
{\bf 2.1.~Proofs of Theorems~\hl{T1}{1}, \hl{T2}{2}, \hl{T4}{4} and \hl{T5}{5}.} We first consider the case of Theorems~\hl{T1}{1} and \hl{T2}{2}. By Proposition~\hl{P1.3}{1.3}, the two conditions in (\hl{1}{1}) are equivalent respectively to
\htt{2.1.1}{}
$$\aligned&\Gr_{\Vt}^{\be}\OO_U\eq0\q\h{if}\q\be\slt k{+}1,\\&\Gr_{\Vt}^{\be}\OO_U\eq0\q\h{if}\q\be\less k{+}1.\endaligned
\leqno(2.1.1)$$
\par\nin Here $U\sst Y$ is any open subvariety on which a local defining function $f$ of $X\sst Y$ is defined.
\sk
By Corollary~\hl{C1.3}{1.3}, the two conditions in (\hl{2.1.1}{2.1.1}) are further equivalent respectively to
\htt{2.1.2}{}
$$\aligned\Gr^F_p\Mtf\ual\eq0&\q\h{if}\q p\less k,\,\al\in(0,1)\,\,\,\h{or}\,\,\,p\slt k,\,\al\eq1,\\
\Gr^F_p\Mtf\ual\eq0&\q\h{if}\q p\less k,\,\al\in(0,1],
\endaligned
\leqno(2.1.2)$$
\par\nin for any local defining function $f$ (setting $\be\eq\al\pl p$).
\sk
By Corollary~\hl{C1.2}{1.2}, these are equivalent respectively to
\htt{2.1.3}{}
$$\aligned\Gr^F_p\VX^{\{\la\}}\eq0&\q\h{if}\q p\less k{+}1,\,\la\nes1\,\,\,\,\h{or}\,\,\,\,p\less k,\,\la\eq1,\\
\Gr^F_p\VX^{\{\la\}}\eq0&\q\h{if}\q p\less k{+}1,\,|\la|\eq1.
\endaligned
\leqno(2.1.3)$$
\par\nin \sk
Recall that the $i$\1th component of the filtered de Rham complex is defined by taking the tensor product with $\Om_Y^{i+d_Y}$, where the Hodge filtration $F$ is shifted depending on the degree; more precisely, it is shifted by ${-}i{-}d_Y$ for $i\ins[-d_Y,0]$, see also \cite[(1.2.3)]{JKSY3}. One can easily verify {\it inductively\1} that the above two conditions are equivalent respectively to
\htt{2.1.4}{}
$$\aligned\Gr_F^p\DR_Y\bl(\VX^{\{\la\}}\br)\eq0&\q\h{if}\q p\gess d_X{-}k,\,\la\nes1\,\,\,\,\h{or}\,\,\,\,p\gess d_Y{-}k,\,\la\eq1,\\\Gr_F^p\DR_Y\bl(\VX^{\{\la\}}\br)\eq0&\q\h{if}\q p\gess d_X{-}k,\,|\la|\eq1.
\endaligned
\leqno(2.1.4)$$
\par\nin Here the isomorphisms are considered in $D^b_{\rm coh}(\OO_Y)$, and $F^p\eq F_{-p}$ ($p\ins\Z$).
\sk
By Corollary~\hl{C1.1}{1.1}, the two conditions in (\hl{2.1.4}{2.1.4}) are respectively equivalent to
\htt{2.1.5}{}
$$\aligned\Gr_F^p\DR_Y\bl(\VX^{\{\la\}}\br)\eq0&\q\h{if}\q p\less k,\,|\la|\eq1,\\\Gr_F^p\DR_Y\bl(\VX^{\{\la\}}\br)\eq0&\q\h{if}\q p\less k,\,\la\nes1\,\,\,\,\h{or}\,\,\,\,p\less k{+}1,\,\la\eq1.
\endaligned
\leqno(2.1.5)$$
\par\nin (Note that these conditions are local on $Y$.) We then get the proof of the first part of Theorems~\hl{T1}{1} and \hl{T2}{2}, since (\hl{6}{6}) is equivalent to its unnormalized version by (\hl{2}{2}).
\sk
To show the converse we have to use the topological filtration on the Grothendieck group of coherent sheaves which is explained in \hl{2.2}{2.2} below. The argument is by increasing induction on $k$ as is explained in the end of the introduction, using Proposition~\hl{P1}{1} and also (\hl{2.2.2}{2.2.2}) below.
\sk
Assume (\hl{5}{5}) in Theorem~\hl{T1}{1}, which is equivalent to (\hl{4}{4}) by (\hl{2}{2}). By the inductive hypothesis, we see that $\Gr_F^k\DR_Y\bl(\VX^{\{\la\}}\br)$ is {\it locally\1} isomorphic to the dual of a coherent sheaf. (Note that the duality isomorphisms are {\it not\1} given globally on $Y$.) We can then apply (\hl{2.2.2}{2.2.2}) below, and the assertion follows. The argument is similar for Theorem~\hl{T2}{2}, where (\hl{6}{6}) is equivalent to its unnormalized version by (\hl{2}{2}). As for (\hl{7}{7}), we can employ Proposition~\hl{P1}{1}, where duality and (\hl{2.1.5}{2.1.5}) are {\it not\1} needed, since (\hl{2.1.4}{2.1.4}) is sufficient. The argument is essentially the same for Theorems~\hl{T4}{4} and \hl{T5}{5} using the self-duality (\hl{1.1.5}{1.1.5}) and the {\it inductive argument\1} before (\hl{2.1.4}{2.1.4}) together with (\hl{2.2.2}{2.2.2}) below. This finishes the proofs of Theorems~\hl{T1}{1}, \hl{T2}{2}, \hl{T4}{4} and \hl{T5}{5}.
\par\htt{2.2}{}\msn
{\bf 2.2.~Topological filtration} (see also \cite[1.3 (or 1.6 in the preprint version)]{MSS3}). Let $K_0(X)$ be the Grothendieck group of coherent sheaves on a complex algebraic variety $X$. It has the topological filtration, which is denoted by $G$ in this paper, and such that $G_kK_0(X)$ is generated by the classes of coherent sheaves $\Fc$ with $\dim{\rm supp}\,\Fc\less k$, see for instance \cite[Examples~1.6.5 and 15.1.5]{Fu} and \cite{SGA6}. The Todd class transformation $td_*$ induces the isomorphisms
\htt{2.2.1}{}
$$\aligned td_*:K_0(X)_{\Q}&\simto\mopl_k\,{\rm CH}_k(X)_{\Q},\\
\Gr^G_k\1td_*:\Gr_k^GK_0(X)_{\Q}&\simto{\rm CH}_k(X)_{\Q}.\endaligned
\leqno(2.2.1)$$
\par\nin see for instance \cite[Corollary 18.3.2]{Fu}. Here we define the increasing filtration $G_p$ on $\mopl_k\,{\rm CH}_k(X)_{\Q}$ by taking the direct sum over $k\less p$.
\sk
Assume $X\eq\PP^n$. For any coherent sheaf $\Fc$ on $\PP^n$ with $k\eq\dim{\rm supp}\,\Fc$, we have the {\it positivity}\,:
\htt{2.2.2}{}
$$\Gr_k^G[\Fc]\eq\msum_{i=1}^r\,m_i\deg Z_i\sgt0\,\,\,\,\,\h{in}\,\,\,\,\Gr_k^GK_0(\PP^n)_{\Q}\eq{\rm CH}_k(\PP^n)_{\Q}\eq\Q,
\leqno(2.2.2)$$
\par\nin with $Z_i$ highest-dimensional irreducible components of ${\rm supp}\,\Fc$ and $m_i\ins\Z_{>0}$ the multiplicity of $\Fc$ at the generic point of $Z_i$ for $i\ins[1,r]$. Recall that $\deg Z_i$ is defined as the intersection number of $Z_i$ with a sufficiently general linear subspace of the complementary dimension.
\par\htt{2.3}{}\msn
{\bf 2.3.~A variant of Hironaka's example.} Let $\pi\,{:}\,\Yt\tos Y\defs\PP^3$ be the blow-up along the union of three lines $C_i$ such that $C_i\cap C_{i+1}\eq\{P_i\}$ with $P_i\nes P_{i+1}$ ($i\ins\Fb_3$), where $\Fb_3\defs\Z/3\Z$. Around $P_i$ the blow-up is taken first along $C_i$ and then along the {\it proper transform\1} of $C_{i+1}$ ($i\ins\Fb_3$). Set
$$E_i\defs\pi^{-1}(C_i),\q D_i\defs\pi^{-1}(P_i),\q A_i\defs D_i\cap E_{i+1},\q B_i\defs\overline{D_i\stm A_i}.$$
\par\nin By the same argument as in Hironaka's example \cite{Hi}, \cite[Appendix B, Example 3.4.1]{Ha}, we get that
\htt{2.3.1}{}
$$\msum_{i\in\Fb_3}\,[B_i]=0\q\h{in}\,\,H_2(E,\Q)\q\h{with}\q E\defs\mcup_i\,E_i.
\leqno(2.3.1)$$
\par\nin Indeed, $[A_0]=[A_1]\pl[B_1]=[A_2]\pl[B_2]\pl[B_1]=[A_0]\pl[B_0]\pl[B_2]\pl[B_1]$.
\sk
There are hyperplanes $H',H_i\subset\PP^3$ such that $H'\,{\cap}\,H_i\eq C_i$ ($i\ins\Fb_3$). Note that their union is a divisor with normal crossings. Take $\si_i,\si'\ins\Gamma\bl(\PP^3,\OO_{\PP^3}(1)\br)$ whose zero-loci are $H_i,H'$ respectively $(i\ins\Fb_3$). These can be identified with projective coordinates $z_1,\dots,z_4$ of $\PP^3$.
\sk
Let $q$ be an integer at least 2. Let $\Xt_i\sst\Yt$ be the proper transform of
$$X_i\defs\{\eta_i\defs\si_i^q\mi\si'\1{}^q\eq0\}\sst\PP^3\q(i\ins\Fb_3).$$
\par\nin We see that
\htt{2.3.2}{}
$$\h{the union $\mcup_{i\in\Fb_3}\,X_i$ has only normal crossings outside $\mcup_i\,C_i$.}
\leqno(2.3.2)$$
\par\nin This is shown for its affine cone in $\C^4$ defined by the product of $h_i\defs z_i^q\mi z_4^q$ ($i\ins[1,3]$). We use the logarithmic functions $\log z_i$ ($i\ins[1,4]$) as analytic local coordinates of $(\C^*)^4$. It is easy to verify the assertion on a {\it sufficiently small\1} neighborhood of each point of $\mcup_{i=1}^4\{z_i\eq0\}$.
\sk
Let $\Xt\sst\Yt$ be the union of the $\Xt_i$ and $E_i$ ($i\ins\Fb_3$). We can prove that
\htt{2.3.3}{}
$$\h{the divisor $\Xt$ has only normal crossings outside $\mcup_i\,B_i$.}
\leqno(2.3.3)$$
\par\nin Indeed, {\it fix\1} $i\ins\Fb_3$, and set
$$x\eq\si_i/\si_{i-1},\q y\eq\si_{i+1}/\si_{i-1},\q z\eq\si'/\si_{i-1},\q g_j\eq\eta_j/\si_{i-1}^q,$$
\par\nin on the affine space $\PP^3\stm H_{i-1}$ containing $P_i$. The polynomials $g_{i-1}$, $g_i$, $g_{i+1}$ are expressed respectively as
$$1\mi z^q,\q x^q\mi z^q,\q y^q\mi z^q.$$
\par\nin \sk
Taking the blow-up along $C_i\stm H_{i-1}$ which is defined by $x,z$ in $\PP^3\stm H_{i-1}\eq\C^3$, we have an affine chart such that the pullbacks of $x,y,z$ are $x'$, $y'$, $x'z'$ respectively with exceptional divisor locally defined by $x'$. Here the polynomials $g_{i-1}$, $g_i$, $g_{i+1}$ become respectively
$$1\mi x'{}^qz'{}^q,\q x'{}^q(1\mi z'{}^q),\q y'{}^q\mi x'{}^qz'{}^q.$$
\par\nin Note that this affine chart meets the proper transform of $C_{i+1}$.
\sk
Taking the blow-up along the latter which is defined by $y,z$ in $\C^3$, we have an affine chart such that the pullbacks of $x',y',z'$ are $x''$, $y''$, $y''z''$ respectively with $E_{i+1}$ locally defined by $y''$. Here $g_{i-1}$, $g_i$, $g_{i+1}$ become respectively
$$1\mi x''{}^qy''{}^qz''{}^q,\q x''{}^q(1\mi y''{}^qz''{}^q),\q y''{}^q(1\mi x''{}^qz''{}^q).$$
\par\nin We have another affine chart such that the pullbacks of $x',y',z'$ are $x'''$, $y'''z'''$, $z'''$ respectively with $E_{i+1}$ locally defined by $z'''$. Here $g_{i-1}$, $g_i$, $g_{i+1}$ become respectively
$$1\mi x'''{}^qz'''{}^q,\q x'''{}^q(1\mi z'''{}^q),\q z'''{}^q(y'''{}^q\mi x'''{}^q).$$
\par\nin Combining these with (\hl{2.3.2}{2.3.2}), we can  verify that the assertion (\hl{2.3.3}{2.3.3}) holds. Here we use logarithmic functions on the open subset $(\C^*)^3$ as analytic local coordinates.
\sk
The above calculation moreover shows that the transversal slice to $B_i\sst \Xt$ at a general point is locally defined by $u(u^q\mi v^q)\eq0$ with $u,v$ analytic local coordinates. So the minimal exponent is $2/(q{+}1)<1$ by the assumption on $q$. Combined with (\hl{2.3.1}{2.3.1}), (\hl{2.3.3}{2.3.3}) and using the action of the cyclic group $\Z/3\Z$ on $\Xt$, this gives a counterexample to Theorem~\hl{T5}{5} in the non-projective compact case.
\par\htt{R2.3}{}\msn
{\bf Remark~2.3.} It seems very difficult to construct a counterexample to Theorem~\hl{T5}{5} with $H_{2\dX}(X,\Q)$ replaced by $H_{2\dX}(\Sing\,X,\Q)$. Indeed, if one takes a blow-up of a fourfold with one-dimensional center modifying the above construction, the minimal exponent locus is {\it not\1} one-dimensional as far as tried.
\msn

\sk
{\smaller\smaller
L.G. Maxim: Department of Mathematics, University of Wisconsin-Madison, 480 Lincoln Drive, Madison, WI 53706-1388, USA

{\it Email address}\,: maxim@math.wisc.edu

M. Saito: RIMS Kyoto University, Kyoto 606-8502 Japan

{\it Email address}\,: msaito@kurims.kyoto-u.ac.jp

R. Yang: Max-Planck-Institut f\"ur Mathematik, Vivatsgasse 7, 53111 Bonn, Germany

{\it Email address}\,: ruijie.yang@hu-berlin.de, ryang@mpim-bonn.mpg.de}
\end{document}